\documentclass[12pt]{article}
\usepackage{amsmath,amsfonts,latexsym,amsthm,amssymb,mathabx}
\newcommand{\Rn}{\mathbb R^n}
\newcommand{\Qp}{\mathbb Q_p}
\newcommand{\Qpn}{\mathbb Q_p^n}
\newcommand{\PP}{\mathbf P}

\DeclareMathOperator{\R}{Re}

\newcommand{\A}{\widetilde{a}}
\newcommand{\G}{\mathfrak{G}}
\newcommand{\GG}{\widehat{\mathfrak {G}}}

\newcommand{\D}{\mathcal D}
\newcommand{\LL}{\mathcal L}
\newcommand{\PF}{\mathfrak{P}}

\begin{document}
\newtheorem{prop}{Proposition}
\newtheorem{lem}{Lemma}
\newtheorem{teo}{Theorem}
\pagestyle{plain}
\title{ Ground States for Nonlocal Schr\"odinger Type Operators on Locally Compact Abelian Groups}

\author{ \textbf{Anatoly N.
Kochubei}\\
Institute of Mathematics,\\
National Academy of Sciences of Ukraine, \\
Tereshchenkivska 3, \\
Kyiv, 01004 Ukraine\\
Email: kochubei@imath.kiev.ua \and \textbf{Yuri Kondratiev}\\
Department of Mathematics, University of Bielefeld, \\
D-33615 Bielefeld, Germany,\\
Email: kondrat@math.uni-bielefeld.de}

\numberwithin{equation}{section}
\date{}
\maketitle

\vspace*{3cm}
\begin{abstract}
We find classes of nonlocal operators of Schr\"odinger type on a locally compact noncompact Abelian group $\G$, for which there exists a ground state. In particular, such a result is obtained for the case where the principal part of our operator generates a recurrent random walk. Explicit conditions for the existence of a ground state are obtained for the case $\G =\Qpn$ where $\Qp$ is the field of $p$-adic numbers.
\end{abstract}
\vspace{2cm}
{\bf Key words: }\ ground state; locally compact Abelian group; operators of Schr\"odinger type; recurrent random walk; field of $p$-adic numbers

\medskip
{\bf MSC 2010}. Primary: 47G10. Secondary: 45P05; 11S80, 43A70.

\section{Introduction}

In a number of applications (see the references in \cite{KMV}) we encounter nonlocal Schr\"odinger operators
\begin{equation}
\label{1.1}
(\LL u)(x)=-m(x)u(x)+\int\limits_{\Rn}a(x-y)u(y)\,dy
\end{equation}
where $a(x)\ge 0$ is an even bounded continuous function, such that $\int\limits_{\Rn}a(x)\,dx=1$, and typically $m\in C_b(\Rn )$, $0\le m(x)\le 1$, $m(x)\to 1$ as $|x|\to \infty$. The recent papers \cite{KMPZ1,KMPZ2,KMV} contain a rather well-developed spectral theory of such operators, in many of its features different from classical counterparts. The reason of these differences is the fact that the convolution operator in (\ref{1.1}) is bounded in $L^2(\Rn)$, and the ``potential'' is not subordinated to it.

On the other hand, the above operators make sense in a much more general setting, where the convolution is defined on a locally compact noncompact Abelian group $\G$. Here, in addition to the case where $\G =\Rn$, we have another important example, in which $\G$ is the field $\Qp$ of $p$-adic numbers. The spectral theory available in this framework (see the books \cite{AKS,KKZ,K2001,VVZ}) deals mostly with perturbations of Vladimirov's fractional differentiation operator, that is does not cover the case of a bounded convolution operator.

In this paper, we concentrate, in the group situation, on the problem of existence of a ground state, the positive eigenfunction corresponding to the maximal positive eigenvalue of $\LL$. We show that the approach of \cite{KMPZ1} based on the use of the Krein-Rutman theorem \cite{KR} remains valid in the general group situation. In addition, we derive an existence result from a recurrence property of the random walk corresponding to the density $a$ \cite{PS}. For the $p$-adic case, we find explicit conditions for the existence of a ground state, parallel to some results from \cite{KMPZ1} but possessing certain special features.

\medskip
\section{Spectral properties. Applications of the Krein-Rutman theorem}

let $\G$ be a second countable locally compact noncompact Abelian group, $\GG$ be its dual group of continuous characters. The action of a character $z\in \GG$ on an element $x\in \G$ will be denoted $\langle x,z\rangle$. We write the group operation additively. Let $\mu$ and $\widehat{\mu}$ be the Haar measures on $\G$ and $\GG$ respectively, normalized in such a way that the Plancherel identity holds, that is the Fourier transform
$$
\widetilde{f}(z)=\int\limits_{\G}\langle x,z\rangle f(x)\mu (dx)
$$
defines an isometric isomorphism from $L^2(\G )$ onto $L^2(\GG )$. The inverse operator is given by the expression
$$
f(x)=\int\limits_{\GG}\overline{\langle x,z\rangle} \widetilde{f}(z)\widehat{\mu} (dz).
$$
See \cite{HR,Ru,We} for information regarding harmonic analysis on topological groups.

Following (\ref{1.1}) we consider the operator
\begin{equation}
\label{2.1}
(\LL u)(x)=-m(x)u(x)+\int\limits_{\G}a(x-y)u(y)\,\mu(dy)
\end{equation}
where $a(x)$ is a nonnegative even function on $\G$,
$$
\int\limits_{\G}a(x)\mu(dx)=1,
$$
and the set of zeroes $N(a)=\{ x\in \G:\ a(x)=0\}$ is empty or compact, $m$ belongs to the Banach space $C_b(\G)$ of bounded continuous (in general, complex-valued) functions on $\G$, $0\le m(x)\le 1$, $m(x)\to 1$, as $x\to \infty$ (the notation $x\to \infty$ means the convergence by the filter generated by complements to compact subsets of $\G$). It is often convenient to rewrite (\ref{2.1}) as
$$
\LL u=\LL_0u+V(x)u
$$
where
$$
(\LL_0u)(x)= \int\limits_{\G}a(x-y)[u(y)-u(x)]\,\mu(dy),\quad V(x)=1-m(x),
$$
so that $0\le V(x)\le 1$. We assume that $V\in C_b(\G)$ and $V(x)\to 0$, as $x\to \infty$.

Below we study spectral properties of $\LL$ both in $L^2(\G)$ and $C_b(\G)$.

\medskip
\begin{lem}
Both in $L^2(\G)$ and $C_b(\G)$, the operator $\LL_0$ is bounded and dissipative, that is
\begin{equation}
\label{2.2}
\| \lambda f-\LL_0f\| \ge \lambda \|f\|\ \text{for any $\lambda >0$}.
\end{equation}
\end{lem}

\medskip
{\it Proof}. The case of $L^2(\G)$ follows immediately from the representation $\LL_0f=a*f-f$, properties of the Fourier transform and the fact that in the case of a Hilbert space the condition (\ref{2.2}) is equivalent to the inequality
$$
\R (\LL_0x,x)\le 0 \ \text{for all $x\in L^2(\G)$}.
$$

In the case of $C_b(\G)$, the operator $\LL_0$ is a generator of a strongly continuous contraction semigroup
\begin{equation}
\label{2.pos}
e^{t\LL_0}=e^{-t}\sum\limits_{k=0}^\infty t^k\frac{a^{*k}}{k!},
\end{equation}
which implies (\ref{2.2}). $\qquad \blacksquare$

\medskip
Note that by (\ref{2.pos}), due to the compactness of $N(a)$, $e^{t\LL_0}$ is a semigroup of convolution operators with strictly positive kernels. The semigroup $e^{t\LL}$ can be given by the Feynman-Kac formula valid for general Markov processes; in the generality sufficient for our situation, that is proved in \cite{Ge}. By this formula, $e^{t\LL}$ has, as well as $e^{t\LL_0}$, the positivity improving property (\cite{RS4}, Section XIII.12).

It follows from Lemma 1 that the half-plane
$$
\D =\{ \lambda \in \mathbb C:\ \R \lambda >0\}
$$
belongs to the resolvent set of the operator $\LL_0$. Now we turn to the operator $\LL$.

\medskip
\begin{lem}
The operator $\LL$ has only discrete spectrum in the half-plane $\D$.
\end{lem}

\medskip
{\it Proof}. For any $\lambda \in \D$ we have
\begin{equation}
\label{2.*}
\lambda -\LL_0-V=(\lambda -\LL_0)(I-(\lambda -\LL_0)^{-1}V),
\end{equation}
so that
$$
(\lambda -\LL_0-V)^{-1}=(I-(\lambda -\LL_0)^{-1}V)^{-1}(\lambda -\LL_0)^{-1}.
$$
Writing the Neumann series for $(\lambda -\LL_0)^{-1}$ we find that
\begin{equation}
\label{2.3}
(\lambda -\LL_0)^{-1}=\frac1{\lambda +1}+\frac1{\lambda +1}A_\lambda
\end{equation}
where
\begin{equation}
\label{2.4}
A_\lambda =(\lambda +1)(\lambda -\LL_0)^{-1}-1=\sum\limits_{l=1}^\infty \frac{a^{*l}}{(\lambda +1)^l}.
\end{equation}

$A_\lambda$, $\lambda \in \D$, is a bounded convolution operator with an integral kernel
\begin{equation}
\label{2.5}
G_\lambda (x-y)=\int\limits_{\GG}\overline{\langle x-y,z\rangle}\frac{\A (z)}{\lambda +1-\A (z)}\widehat{\mu}(dz),
\end{equation}
with
\begin{equation}
\label{2.6}
\widetilde{G_\lambda}(z)= \frac{\A (z)}{\lambda +1-\A (z)}\in L^2(\GG ),\quad \lambda \in \D.
\end{equation}
It follows from (\ref{2.4}) that $G_\lambda \in L^1(\G)\cap C_b(\G)$.

Denote by $W_\lambda$ the operator of multiplication by the function
$$
W_\lambda (x)=1-\frac{V(x)}{\lambda +1},\quad \lambda \in \D.
$$
It is a bounded operator with a bounded inverse. It follows from (\ref{2.3}) that
$$
I-(\lambda -\LL_0)^{-1}V=W_\lambda -\frac1{\lambda +1}A_\lambda V,
$$
so that
\begin{equation}
\label{2.7}
(\lambda -\LL_0-V)^{-1}=(I-Q_\lambda)^{-1}((\lambda -\LL_0)W_\lambda)^{-1}
\end{equation}
where
\begin{equation}
\label{2.8}
Q_\lambda =\frac1{\lambda +1}W_\lambda^{-1}A_\lambda V.
\end{equation}

Since $\G$ is locally compact and second countable, there exists its countable covering by open subsets with compact closures. Taking the subordinate partition of unity, we can form a monotone sequence of continuous nonnegative functions $\{ h_r\}$ with compact supports converging pointwise to 1. Let $A_\lambda^{(r)}$ be the convolution operator with the convolution kernel $G_\lambda^{(r)}(x)=h_r(x)G_\lambda (x)$, $V^{(r)}$ be the operator of multiplication by $V^{(r)}(x)=h_r(x)V(x)$. Then $A_\lambda^{(r)}$ is a compact operator, and $A_\lambda^{(r)}V^{(r)}\to A_\lambda V$ by the norm of operators. Therefore the operator $Q_\lambda$ is compact.

Using the analytic Fredholm theorem (see \cite{RS1}, Theorem VI.14, and \cite{St}), we find that $(I-Q_\lambda)^{-1}$ is meromorphic on $\D$. Since $((\lambda -\LL_0)W_\lambda)^{-1}$ is a bounded operator, we conclude that the operator $\LL$ has only a discrete spectrum in $\D$. $\qquad \blacksquare$

\medskip
It follows from (\ref{2.4}), (\ref{2.5}) and (\ref{2.8}) that
\begin{description}
\item [(i)] $A_\lambda$ is a positivity improving operator for each $\lambda >0$, since $G_\lambda (x-y)>0$, $x,y\in \G$.
\item [(ii)] $G_\lambda (x-y)$ is monotonically decreasing with respect to $\lambda >0$;
\item [(iii)] $Q_\lambda$, $\lambda >0$, is a positivity improving compact integral operator on $C_b(\G)$ with the kernel
$$
Q_\lambda (x,y)=\frac{G_\lambda (x-y)V(y)}{\lambda +1-V(x)}.
$$
\end{description}

\medskip
Let us study the spectral radius $r(Q_\lambda)$ of the operator $Q_\lambda$, $\lambda >0$. Let $\mathcal E$ be one of the spaces $C_b(\G)$ or $L^2(\G)$.

\medskip
\begin{lem}[\cite{KMPZ1}]
The spectral radius $r(Q_\lambda)$ is continuous and monotonically decreasing with respect to $\lambda >0$, and $r(Q_\lambda)\to 0$ for $\lambda \to +\infty$. If there exists a function $\varphi \in \mathcal E$, $\varphi \ge 0$, $\| \varphi\|=1$, such that $Q_\lambda \varphi (x)\ge c_0\varphi (x)$, then $r(Q_\lambda)\ge c_0$.
\end{lem}

\medskip
By (\ref{2.*}), the equation on the eigenfunction $\psi$,
\begin{equation}
\label{2.9}
(\LL_0+V-\lambda )\psi =0,\quad \lambda >0,
\end{equation}
is equivalent to the equality
\begin{equation}
\label{2.10}
Q_\lambda \psi =\psi.
\end{equation}
It follows from Lemma 3 that if
\begin{equation}
\label{2.11}
\lim\limits_{\lambda \to +0}r(Q_\lambda )>1,
\end{equation}
then there exists such $\lambda >0$ that
$$
r(Q_\lambda )=1, \quad \text{and $r(Q_{\lambda'} )<1$ for $\lambda'>\lambda$}.
$$

By the Krein-Rutman theorem (\cite{KR}, Theorem 6.1), in this case 1 is the maximal positive eigenvalue of $Q_\lambda$ with a positive eigenfunction $\psi_\lambda >0$. Correspondingly, $\lambda$ is the maximal positive eigenvalue of the operator $\LL$, and $\psi_\lambda (x) >0$ is the ground state of $\LL$. The uniqueness of the ground state follows from the positivity improving property of the semigroup $e^{t\LL}$; see Theorem XIII.44 in \cite{RS4} or Theorem 2.11 in \cite{Kr}.

The following properties of the ground state, provided it exists, are proved just as in \cite{KMPZ1}. The only difference is that instead of the standard dominated convergence theorem, we use its generalization to the case of convergence with respect to a filter with a countable basis (\cite{Bo}, Chapter 4, \S 3.7). In our case, this countable basis can be constructed as follows: take the countable increasing sequence of compact sets appearing in the construction of a partition of unity (see e.g. Theorem 5.3 in \cite{La}), and take the sequence of their complements in $\G$ as the basis of a filter.

\medskip
\begin{lem}
\begin{description}
\item [(i)] If the ground state $\psi_\lambda$ belongs to $C_b(\G)$, then $\psi_\lambda (x)\to 0$, as $x\to \infty$.
\item [(ii)] If the ground state $\psi_\lambda\in L^2(\G)$, then $\psi_\lambda \in C_b(\G)\cap L^2(\G)$ and $\psi_\lambda (x)\to 0$, as $x\to \infty$.
\item [(iii)] If the ground state $\psi_\lambda \in C_b(\G)$ and $V\in L^2(\G)$, then $\psi_\lambda\in L^2(\G)$.
\end{description}
\end{lem}

\medskip
{\it Remark}. Let $\mathcal E_1$ be one of the spaces $C_b(\G),L^2(\G)$, and $\mathcal E_2$ be the other of these spaces. As we see in Lemma 4, it is a typical situation that a ground state in $\mathcal E_1$ belongs also to $\mathcal E_2$. In this case, is it a ground state of the operator in $\mathcal E_2$?

The answer is positive, because in our setting with bounded operators, being an eigenfunction is a pointwise property. By the Krein-Rutman theorem (\cite{KR}, Theorem 6.1), the operator in $\mathcal E_2$, having a positive eigenvalue, has also a ground state. In addition, for positivity improving operators, a normalized positive eigenfunction is unique (\cite{Kr}, Theorem 2.11). Therefore the given ground state for $\mathcal E_1$ is also the ground state for $\mathcal E_2$.

\medskip
\section{Random walk}

Let us consider the probability measure $\mu_a(dx)=a(x)\mu (dx)$ on the Borel $\sigma$-algebra of the group $\G$. If $\xi_1,\xi_2,\ldots ,\xi_m,\ldots$ are independent $\G$-valued random variables, each having $\mu_a$ for its probability law, then the random walk with an initial point $S_0$ is the Markov chain $S_m=S_0+\xi_1+\cdots +\xi_m$. Below we connect some properties of the random walk found in \cite{PS} with the existence of the ground state of the operator $\LL$.

In addition to the classical example of $\Rn$ considered in \cite{KMPZ1}, we have the second main example, in which $\G$ is the additive group of the field $\Qp$ of $p$-adic numbers. Harmonic analysis on $\Qp$ and other local fields is a well-developed branch of contemporary mathematical analysis; see \cite{K2001,Ta,VVZ}.

Recall that $\Qp$ is a completion of the field $\mathbb Q$ of rational numbers with respect to the absolute value $|x|_p$ defined
by setting $|0|_p=0$,
$$
|x|_p=p^{-\nu }\ \mbox{if }x=p^\nu \frac{m}n,
$$
where $\nu ,m,n\in \mathbb Z$, and $m,n$ are prime to $p$. $\Qp$ is a locally compact topological field.

Note that by Ostrowski's theorem there are no absolute values on $\mathbb Q$, which are not equivalent to the ``Euclidean'' one,
or one of $|\cdot |_p$.

The absolute value $|x|_p$, $x\in \mathbb Q_p$, has the following properties:
\begin{gather*}
|x|_p=0\ \mbox{if and only if }x=0;\\
|xy|_p=|x|_p\cdot |y|_p;\\
|x+y|_p\le \max (|x|_p,|y|_p).
\end{gather*}

The latter property called the ultra-metric inequality (or the non-Archi\-me\-dean property) implies the total disconnectedness of $\Qp$ in the topology
determined by the metric $|x-y|_p$, as well as many unusual geometric properties. Note also the following consequence of the ultra-metric inequality:
\begin{equation*}
|x+y|_p=\max (|x|_p,|y|_p)\quad \mbox{if }|x|_p\ne |y|_p.
\end{equation*}

The absolute value $|x|_p$ takes the discrete set of non-zero
values $p^N$, $N\in \mathbb Z$. If $|x|_p=p^N$, then $x$ admits a
(unique) canonical representation
\begin{equation}
\label{3.1}
x=p^{-N}\left( x_0+x_1p+x_2p^2+\cdots \right) ,
\end{equation}
where $x_0,x_1,x_2,\ldots \in \{ 0,1,\ldots ,p-1\}$, $x_0\ne 0$.
The series converges in the topology of $\mathbb Q_p$. For
example,
$$
-1=(p-1)+(p-1)p+(p-1)p^2+\cdots ,\quad |-1|_p=1.
$$
We denote $\mathbb Z_p=\{ x\in \Qp:\ |x|_p\le 1\}$. $\mathbb Z_p$, as well as all balls in $\Qp$, is simultaneously open and closed.

Proceeding from the canonical representation (\ref{3.1}) of an element $x\in
\mathbb Q_p$, we define the fractional part of $x$ as the rational number
$$
\{ x\}_p=\begin{cases}
0,& \text{if $N\le 0$ or $x=0$};\\
p^{-N}\left( x_0+x_1p+\cdots +x_{N-1}p^{N-1}\right) ,& \text{if
$N>0$}.\end{cases}
$$
The function $\chi (x)=\exp (2\pi i\{ x\}_p)$ is an additive
character of the field $\mathbb Q_p$, that is a character of its additive group. It is clear
that $\chi (x)=1$ if $|x|_p\le 1$. Denote by $dx$ the Haar measure on the
additive group of $\Qp $ normalized by the equality $\int_{\mathbb Z_p}dx=1$.
The above additive group is self-dual -- every continuous character can be written as $\langle x,z\rangle =\chi (zx)$, $x\in \Qp$, where $z\in \Qp$ is a unique element defining the character.

Returning to random walk on a general group $\G$ we describe the assumption from \cite{PS}, under which the result we need was obtained.

\medskip
$(\G_0)$: \ The minimal closed subgroup of $\G$ generated by $\operatorname{supp} \mu_a$ coincides with $\G$.

\medskip
If $\A (z)=1$ where $z\ne 0$ (note that 0 is the unit character; this additive notation corresponds, for example, to the cases of $\mathbb R$ or $\Qp$), then
$$
1-\A (z)=\int\limits_{\G}(1-\langle x,z\rangle )a(x)\mu (dx),
$$
so that $\langle x,z\rangle=1$ for $x\in \operatorname{supp} \mu_a$, thus for all $x$ from the subgroup generated by $\operatorname{supp} \mu_a$, so that due to $(\G_0)$, $\langle x,z\rangle=1$ for all $x\in \G$. This contradiction proves that
\begin{equation}
\label{3.2}
\A (z) \ne 1 \ \text{for $z\ne 0$}.
\end{equation}

\medskip
For $\G=\mathbb R$ (or $\Rn$), the measures whose Fourier transforms equal 1 at nonzero points (or slightly more general ones) are called the lattice distributions and are concentrated on lattices (see Section 1.5 in \cite{Ra} or Section XV.1 in \cite{Fe}). This is of course impossible in the case of a continuous density.

For the case of $\Qp$ (or, similarly, for $\Qpn$) suppose that $\operatorname{supp}a\subset \left\{ x\in \Qp :\ |x|_p\le p^N\right\}$, $N\in \mathbb Z$. We have
$$
1-\A (\xi )=\int\limits_{|x|_p\le p^N} [1-\cos (2\pi \{ x\xi \}_p)]a(x)\,dx,
$$
so that $1-\A (\xi )=0$ for $|\xi|\le p^{-N}$. Thus, in the $p$-adic case, the property (\ref{3.2}) is violated, if $a$ has a compact support. Conversely, if $\A (\xi)=1$ for $\xi \ne 0$, then $a(x)=0$ for $|x|_p>|\xi |_p^{-1}$. On the other hand, this case is excluded by the condition $(\G_0)$ -- by the ultrametric inequality, a $p$-adic ball does not generate the additive group of $\Qp$.

A random walk is said to be recurrent, if for some compact neighborhood $M$ of 0,
$$
\sum\limits_{m=1}^\infty \PP \{ S_m\in M|S_0=0\} =\infty.
$$

\medskip
\begin{teo}
Suppose that $\G$ is a noncompact second countable locally compact Abelian group, and $a(x)$, in addition to the assumptions made in Introduction, is such that $(\G_0)$ is satisfied, and the corresponding random walk is recurrent. Then for any $V\not\equiv 0$ satisfying the conditions from Section 2, the ground state of the operator $\LL$ in $C_b(\G)$ exists.
\end{teo}

\medskip
{\it Proof}. It follows from the recurrence of our random walk that
\begin{equation}
\label{3.3}
\int\limits_P\frac1{1-\A (z)}\widehat{\mu} (dz)=\infty
\end{equation}
for any open neighborhood $P$ of the origin in $\GG$ (\cite{PS}, Theorem 5.1).

Let us take a continuous real-valued function $\varphi$ with a compact support, $\|\varphi \|=1$, such that
$$
\widetilde{V\varphi}(0)=\int\limits_{\G}V(x)\varphi (x)\,\mu (dx)>0.
$$

The assumptions regarding $a(x)$ and $V(x)$ imply that $\A ,\widetilde{V\varphi}\in L^2(\G)\cap C_b(\G)$ and consequently, $\A \widetilde{V\varphi}\in L^1(\G)\cap C_b(\G)$. Note also that $\A (z)\to 0$, as $z\to \infty$ (\cite{Ru}, Theorem 1.2.4). Then by (\ref{2.5}), for any $x\in \operatorname{supp}\varphi$,
\begin{multline*}
(Q_\lambda \varphi )(x)=\frac1{\lambda +1-V(x)}\int\limits_{\G}G_\lambda (x-y)V(y)\varphi (y)\,\mu (dy)\\
=\frac1{\lambda +1-V(x)}\int\limits_{\G}V(y)\varphi (y)\,\mu (dy)\int\limits_{\GG}\langle x-y,z\rangle \frac{\A (z)}{\lambda +1-\A (z)}\widehat{\mu}(dz)\\
=\frac1{\lambda +1-V(x)}\int\limits_{\GG}\frac{ \langle -x,z\rangle \A (z)\widetilde{V\varphi}(z)}{\lambda +1-\A (z)}\widehat{\mu}(dz).
\end{multline*}

Since all functions appearing in the above expression of $Q_\lambda \varphi$ as an integral on $\G$ are real-valued, we may write
\begin{equation}
\label{3.4}
(Q_\lambda \varphi )(x)=\frac1{\lambda +1-V(x)}\int\limits_{\GG}\frac{\A (z)\R \left\{  \langle -x,z\rangle \widetilde{V\varphi}(z)\right\}}{\lambda +1-\A (z)}\widehat{\mu}(dz).
\end{equation}

Choose such a neighborhood $\PF$ of the origin in $\GG$ that
\begin{equation}
\label{3.5}
\A (z)\R \left\{  \langle -x,z\rangle \widetilde{V\varphi}(z)\right\} \ge \gamma_0>0,\quad z\in \PF,
\end{equation}
and write the integral in (\ref{3.4}) as the sum of integrals over $\PF$ and $\GG \setminus \PF$. Due to (\ref{3.2}), the second of them has a finite limit, as $\lambda \to 0$. By the monotone convergence theorem, it follows from (\ref{3.3}) and (\ref{3.4}) that the integral over $\PF$ tends to $+\infty$, as $\lambda \to +0$.

The continuous function $U_\lambda (x)=[(Q_\lambda \varphi )(x)]^{-1}$ tends to 0 monotonically, as $\lambda \to 0$ (this monotonicity follows from the monotone dependence of $G_\lambda$ on $\lambda$). By the Dini theorem, this convergence is uniform on $\operatorname{supp}\varphi$. Therefore for any $c_0>0$, there exists such $\lambda >0$ that $(Q_\lambda \varphi (x)\ge c_0\varphi (x)$. By Lemma 3, this means that $\lim\limits_{\lambda \to +0} r(Q_\lambda )=\infty$, which implies the inequality (\ref{2.11}) and the existence of a ground state. $\qquad \blacksquare$

\medskip
\section{The $p$-adic case}

In this section we give explicit conditions for the existence of a ground state in the case where $G=\Qpn$. We consider mostly the case of the ground state from $C_b(\Qpn)$. The existence of a ground state for operators on $L^2(\Qpn)$ can be proved just as it is done for operators on $L^2(\Rn)$ in (\cite{KMPZ1}, Section 4).

Let us introduce some notation for $\Qpn =\underset{\text{$n$ times}}{\underbrace{\Qp \times \cdots \times \Qp}}$. This is a vector space over the field $\Qp$ with the norm
$$
\| x\| =\max\limits_{1\le j\le n}|x_j|_p,\quad x=(x_1,\ldots ,x_n)\in \Qpn.
$$
$B_N$ ($N\in \mathbb Z$) will denote the ball $\left\{ x\in \Qpn :\ \|x\| \le p^N\right\}$. Just as for $\Qp$, the additive group of $\Qpn$ is self-dual. The simplest integration formulas are
$$
\int\limits_{\|x\| \le p^j}d^n x=p^{nj},\quad \int\limits_{\|x\| = p^j}d^n x=(1-p^{-n}) p^{nj},\quad j\in \mathbb Z
$$
(see e.g. \cite{Vl}). The Fourier transform of a function $f:\ \Qpn \to \mathbb C$ is defined as
$$
\widetilde{f}(\xi )=\int\limits_{\Qpn}\chi (x\cdot \xi )f(x)\,d^n x
$$
where $x\cdot \xi =\sum\limits_{j=1}^n x_j\xi_j$. As before, the Fourier transform of an even function, for example of $a$, is real-valued.

The next result gives a simple example of the existence of a ground state.

\medskip
\begin{teo}
Suppose that $V(x)=1$ for $x\in B_N$, $N\in \mathbb Z$. Then the ground state of $\LL$ exists.
\end{teo}

\medskip
{\it Proof}. Let $f_N$ be the indicator of the ball $B_N$. Since $B_N$ is simultaneously open and closed, $f_N\in C_b(\Qpn)$. For any $\lambda \in (0,1)$ and any $x\in B_N$, we get
$$
(Q_\lambda f_N)(x)=\int\limits_{B_N}Q_\lambda (x,y)\,dy\ge \frac1\lambda \int\limits_{B_N}G_\lambda (x-y)\,dy\ge \frac{\operatorname{Vol}(B_N)}{\lambda}\varkappa
$$
where $\varkappa =\min\limits_{x,y\in B_N}G_1(x-y)<\min\limits_{\lambda \in (0,1)}\min\limits_{x,y\in B_N}G_\lambda (x-y)$.

By Lemma 3, $\lim\limits_{\lambda \to +0}r(Q_\lambda)=\infty$, which implies the existence of a ground state. $\qquad \blacksquare$

\medskip
In our next result, we begin with a ground state in $L^2(\Qpn)$, and then use Remark after Lemma 4.

\medskip
\begin{teo}
Assume that for some $\beta \in (0,1)$, there exists such $N\in \mathbb Z$ that
$$
\beta \le V(x)\le 1, \quad x\in B_N.
$$
Then the ground state of the operator $\LL$ exists, if $N=N(\beta)$ is sufficiently large.
\end{teo}

\medskip
{\it Proof}. To prove the existence of a ground state $\psi_\lambda \in L^2(\Qpn)$, it suffices to show that the quadratic form $(\LL f,f)$ is positive for some $f\in L^2(\Qpn)$. Taking the indicator $f_N$, as in the proof of Theorem 2, we have
$$
(\LL f_N,f_N)=(\LL_0 f_N,f_N)+(Vf_N,f_N)
$$
and
\begin{equation}
\label{4.1}
(Vf_N,f_N)\ge \beta \operatorname{Vol}(B_N).
\end{equation}

Next, if $\|x\| \le p^N$, then
$$
(\LL_0 f_N)(x)=-\int\limits_{\|y\| >p^N}a(x-y)\,dy,
$$
so that
$$
-(\LL_0 f_N,f_N)=\int\limits_{\|x\|\le p^N}dx\int\limits_{\|y\| >p^N}a(x-y)\,dy.
$$
The change of variables $y=x-z$ implies (note that $\|z\|=\|y\|$) the equality
$$
-(\LL_0 f_N,f_N)=\int\limits_{\|x\|\le p^N}dx \int\limits_{\|z\| >p^N}a(z)\,dz
$$
and the relation
\begin{equation}
\label{4.2}
\frac1{\operatorname{Vol}(B_N)}(\LL_0 f_N,f_N)\to 0,\quad \text{as $N\to \infty$}.
\end{equation}

Comparing (\ref{4.1}) and (\ref{4.2}) we find that $\LL$ has a positive discrete spectrum, hence a ground state, if $N$ is large enough. $\qquad \blacksquare$

\medskip
By our assumption, $N(a)$ is empty or compact. Hence the function $a$ cannot have a compact support, so that $\A (\xi )\ne 1$ for $\xi \ne 0$. Let us study the behavior of $\A (\xi )$ near the origin.

Denote
$$
A_j=\sup\limits_{\|x\|=p^j}a(x),\quad j\in \mathbb Z.
$$

\medskip
\begin{lem}
If for some $l\ge 0$,
\begin{equation}
\label{4.3}
\sum \limits_{N=l}^\infty \frac1{p^{Nn}\sum\limits_{j=N+1}^\infty p^{jn}A_j}=\infty,
\end{equation}
then
\begin{equation}
\label{4.4}
\int\limits_{\|\xi \|\le p^{-l}}\frac{d^n\xi}{1-\A (\xi )}=\infty.
\end{equation}
\end{lem}

\medskip
{\it Proof}. Since the function $a$ is even, $\A$ is real-valued. We have
\begin{multline*}
1-\A(\xi )=\int\limits_{Qpn}[1-\chi (\xi \cdot x)]a(x)\,d^nx=\sum \limits_{j=-\infty}^\infty \int\limits_{\|x\|=p^j}[1-\R \chi (\xi \cdot x)]a(x)\,d^nx\\
\le \sum \limits_{j=-\infty}^\infty A_j\int\limits_{\|x\|=p^j}[1-\R \chi (\xi \cdot x)]\,d^nx.
\end{multline*}

The integral $\int\limits_{\|x\|=p^j}\R \chi (\xi \cdot x)\,d^nx=\int\limits_{\|x\|=p^j}\chi (\xi \cdot x)\,d^nx$ is given in \cite{Vl} (formula (15.10):
$$
\int\limits_{\|x\|=p^j}\chi (\xi \cdot x)\,d^nx =(1-p^{-n})p^{jn}\cdot \begin{cases}
1, & \text{if $\|\xi \|\le p^{-j}$},\\
0, & \text{otherwise}; \end{cases}
-p^{(j-1)n}\cdot \begin{cases}
1, & \text{if $\|\xi \|=p^{-j+1}$},\\
0, & \text{otherwise}. \end{cases}
$$

Let $\|\xi \|=p^{-N}$, $N\ge 0$. If $N\ge j$, then
$$
\int\limits_{\|x\|=p^j}\chi (\xi \cdot x)\,d^nx=(1-p^{-n})p^{jn}.
$$
If $N=j-1$, then
$$
\int\limits_{\|x\|=p^j}\chi (\xi \cdot x)\,d^nx=-p^{(j-1)n}=-p^{Nn}.
$$
For other values of $j$, the integral equals 0. Therefore
$$
\int\limits_{\|x\|=p^j}[1-\chi (\xi \cdot x)]\,d^nx=\begin{cases}
0, & \text{if $j\le N$},\\
(1-p^{-n})p^{jn}+p^{Nn}=p^{(N+1)n}, & \text{if $j=N+1$},\\
(1-p^{-n})p^{jn}, & \text{if $j\ge N+2$},\end{cases}
$$
and we obtain the inequality
$$
0\le 1-\A (\xi )\le \sum\limits_{j=N+1}^\infty p^{jn}A_j,
$$
hence the inequality
$$
\int\limits_{\|\xi \|\le p^{-l}}\frac{d^n\xi }{1-\A (\xi )}\ge (1-p^{-n})\sum\limits_{N=l}^\infty \frac1{p^{Nn}\sum\limits_{j=N+1}^\infty p^{jn}A_j}.
$$

Now the condition (\ref{4.3}) implies (\ref{4.4}). $\qquad \blacksquare$

\medskip
The next theorem follows from Lemma 5, just as Theorem 1 was deduced from a result obtained in \cite{PS}.

\medskip
\begin{teo}
For any functions $a$ and $V\not\equiv 0$ satisfying the assumptions from Introduction and Section 2, as well as the condition (\ref{4.3}), a ground state from $C_b(\Qpn )$ exists.
\end{teo}

\medskip
{\it Remark}. If the function $a$ is radial, that is $a(x)=\mathcal A (\|x\|)$, then $A_j=\mathcal A(p^j)$ and
$$
\sum\limits_{j=N+1}^\infty p^{jn}A_j=(1-p^{-n})^{-1}\int\limits_{\|x\|\ge p^{N+1}}a(x)\,d^nx.
$$

\bigskip
\section*{Acknowledgments}
The first-named author is grateful to the Bielefeld University for hospitality during his visits to Bielefeld. The work of the first author was also supported in part by Grant 23/16-18 ``Statistical dynamics, generalized Fokker-Planck equations, and their applications in the theory of complex systems'' of the Ministry of Education and Science of Ukraine. The second-named author gratefully acknowledges the financial support by the DFG through CRC 701 "Stochastic Dynamics: Mathematical Theory and Applications" and the European Commission under the project STREVCOMS PIRSES-2013-
612669.


\medskip
\begin{thebibliography}{999}
\bibitem{AKS}
S. Albeverio, A. Yu. Khrennikov and V. M. Shelkovich, {\it Theory of p-Adic Distributions. Linear and Nonlinear Models}. Cambridge University Press, 2010.
\bibitem{Bo}
N. Bourbaki, {\it Elements of Mathematics. Integration I}, Springer, Berlin, 2004.

\bibitem{Fe}
W. Feller, {\it An Introduction to Probability Theory and Its Applications}, Vol. 2, Wiley, New York, 1971.
\bibitem{Ge}
R. K. Getoor, Additive functionals of a Markov process, {\it Pacif. J. Math.} {\bf 7} (1957), 1577--1591.
\bibitem{HR}
E. Hewitt and K. A. Ross, {\it Abstract Harmonic Analysis, Vols. I, II}, Springer, Berlin, 1963, 1970.
\bibitem{KKZ}
A. Yu. Khrennikov, S. V. Kozyrev and W. Z\'u\~niga-Galindo, {\it Ultrametric Pseudo-Differential Equations and Applications}, Cambridge University Press, 2018.
\bibitem{K2001}
A. N. Kochubei, {\it Pseudo-Differential Equations and Stochastics
over Non-Archimedean Fields}, Marcel Dekker, New York, 2001.
\bibitem{KMPZ1}
Yu. Kondratiev, S. Molchanov, S. Pirogov and E. Zhizhina, On ground state of some non local Schrodinger operators, {\it Appl. Anal.} {\bf 96} (2017), 1390--1400.
\bibitem{KMPZ2}
Yu. Kondratiev, S. Molchanov, A. Piatnitski and E. Zhizhina, Resolvent bounds for jump generators. {\it Appl. Anal.} {\bf 97} (2018), 323--336.
\bibitem{KMV}
Yu. Kondratiev, S. Molchanov and B. Vainberg, Spectral analysis of non-local Schrodinger operators, {\it J. Funct. Anal.} {\bf 273} (2017), 1020--1048.
\bibitem{Kr}
M. A. Krasnosel'skii, {\it Positive Solutions of Operator Equations}, Noordhoff, Groningen, 1964.
\bibitem{KR}
M. G. Krein and M. A. Rutman, Linear operators leaving invariant a cone in a Banach space, {\it Amer. Math. Soc. Translation}, no. 26 (1950), 128 pp.
\bibitem{La}
S. Lang, {\it Real and Functional Analysis}, Springer, New York, 1993.
\bibitem{PS}
S. C. Port and C. J. Stone, Potential theory of random walks on Abelian groups, {\it Acta Math.} {\bf 122} (1969), 19--114.
\bibitem{Ra}
B. Ramachandran, {\it Advanced Theory of Characteristic Functions}, Statistical Publishing Society, Calcutta, 1967.
\bibitem{RS1}
M. Reed and B. Simon, {\it Methods of Modern Mathematical Physics, Vol. 1}, Academic Press, New York, 1972.
\bibitem{RS4}
M. Reed and B. Simon, {\it Methods of Modern Mathematical Physics, Vol. 4}, Academic Press, New York, 1978.
\bibitem{Ru}
W. Rudin, {\it Fourier Analysis on Groups}, Interscience, New York, 1962.
\bibitem{St}
S. Steinberg, Meromorphic families of compact operators, {\it Arch. Rat. Mech. Anal.} {\bf 31} (1968), 372--379.
\bibitem{Ta}
M. H. Taibleson, {\it Fourier Analysis on Local Fields}, Princeton University Press, 1975.
\bibitem{VVZ}
V. S. Vladimirov, I. V. Volovich and E. I. Zelenov, {\it $p$-Adic Analysis and
Mathematical Physics}, World Scientific, Singapore, 1994.
\bibitem{Vl}
V. S. Vladimirov, {\it Tables of Integrals of Complex-Valued
Functions of $p$-Adic Arguments}, Steklov Mathematical Institute,
Moscow, 2003 (Russian). English version, ArXiv: math-ph/9911027.
\bibitem{We}
A. Weil, {\it L'Int\'egration dans les Groupes Topologiques et ses Applications}, Hermann, Paris, 1965.
\end{thebibliography}
\end{document}